\documentclass[11pt]{article}

\usepackage{latexsym}
\usepackage{graphicx}
\usepackage{float}
\usepackage{amsmath}
\usepackage{amssymb}
\usepackage{lscape}
\usepackage{mathtools}
\usepackage{amsthm}

\usepackage[utf8]{inputenc}
\usepackage{times}

\begin{document}
\bibliographystyle{plain}
\floatplacement{table}{H}
\theoremstyle{theorem}
\newtheorem{theorem}{Theorem}
\newtheorem{lemma}{Lemma}
\newtheorem{corollary}{Corollary}

\theoremstyle{definition}
\newtheorem{definition}{Definition}
\newtheorem{example}{Example}

\newcommand{\sni}{\sum_{i=1}^{n}}
\newcommand{\snj}{\sum_{j=1}^{n}}
\newcommand{\smj}{\sum_{j=1}^{m}}
\newcommand{\sumjm}{\sum_{j=1}^{m}}
\newcommand{\bdis}{\begin{displaymath}}
\newcommand{\edis}{\end{displaymath}}
\newcommand{\beq}{\begin{equation}}
\newcommand{\eeq}{\end{equation}}
\newcommand{\beqn}{\begin{eqnarray}}
\newcommand{\eeqn}{\end{eqnarray}}
\newcommand{\simleq}{\stackrel{<}{\sim}}
\newcommand{\sep}{\;\;\;\;\;\; ; \;\;\;\;\;\;}
\newcommand{\sqrn}{\sqrt{n}}
\newcommand{\sqrtwo}{\sqrt{2}}

\newcommand{\onehlf}{\frac{1}{2}}
\newcommand{\donehlf}{\dfrac{1}{2}}
\newcommand{\thrhlf}{\frac{3}{2}}
\newcommand{\fivhlf}{\frac{5}{2}}
\newcommand{\onethd}{\frac{1}{3}}
\newcommand{\lb}{\left ( }
\newcommand{\lcb}{\left \{ }
\newcommand{\lsb}{\left [ }
\newcommand{\labs}{\left | }
\newcommand{\rb}{\right ) }
\newcommand{\rcb}{\right \} }
\newcommand{\rsb}{\right ] }
\newcommand{\rabs}{\right | }
\newcommand{\lnm}{\left \| }
\newcommand{\rnm}{\right \| }
\newcommand{\lambdab}{\bar{\lambda}}
%
%
\newcommand{\xj}{x_{j}}
\newcommand{\xjb}{\bar{x}_{j}}
\newcommand{\anmell}{a_{n-\ell}}
\newcommand{\anmellmo}{a_{n-\ell-1}}
\newcommand{\anmkmo}{a_{n-k-1}}
\newcommand{\anmkmell}{a_{n-k-\ell}}
\newcommand{\anmkmellmo}{a_{n-k-\ell-1}}
\newcommand{\azero}{a_{0}}
\newcommand{\aone}{a_{1}}
\newcommand{\atwo}{a_{2}}
\newcommand{\ath}{a_{3}}
\newcommand{\afr}{a_{4}}
\newcommand{\afv}{a_{5}}
\newcommand{\asx}{a_{6}}
\newcommand{\anmo}{a_{n-1}}
\newcommand{\anmt}{a_{n-2}}
\newcommand{\an}{a_{n}}
\newcommand{\aj}{a_{j}}
\newcommand{\ajpo}{a_{j+1}}
\newcommand{\ajpt}{a_{j+2}}
\newcommand{\ajmo}{a_{j-1}}
\newcommand{\ajmt}{a_{j-2}}
\newcommand{\ak}{a_{k}}
\newcommand{\akpo}{a_{k+1}}
\newcommand{\akpt}{a_{k+2}}
\newcommand{\akmo}{a_{k-1}}
\newcommand{\akmt}{a_{k-2}}
\newcommand{\aell}{a_{\ell}}
\newcommand{\bone}{b_{1}}
\newcommand{\btwo}{b_{2}}
\newcommand{\defeq}{\coloneqq}      

\newcommand{\lnorm}{\left \|}
\newcommand{\rnorm}{\right \|}
\newcommand{\lnrm}{\biggl | \biggl |}
\newcommand{\rnrm}{\biggr |\biggr |}

\newcommand{\recipp}{p^{\protect \#}}
\newcommand{\recipP}{P^{\protect \#}}

\newcommand{\gammazero}{\gamma_{0}}
\newcommand{\gammaone}{\gamma_{1}}
\newcommand{\sigmaone}{\sigma_{1}}
\newcommand{\sigmatwo}{\sigma_{2}}

\newcommand{\alfazero}{\alpha_{0}}
\newcommand{\alfaone}{\alpha_{1}}
\newcommand{\alfatwo}{\alpha_{2}}
\newcommand{\alfathr}{\alpha_{3}}
\newcommand{\alfath}{\alpha_{3}}
\newcommand{\alfafr}{\alpha_{4}}
\newcommand{\alfafv}{\alpha_{5}}
\newcommand{\alfasx}{\alpha_{6}}
\newcommand{\alfanmo}{\alpha_{n-1}}
\newcommand{\alfanmt}{\alpha_{n-2}}
\newcommand{\alfan}{\alpha_{n}}
\newcommand{\alfaj}{\alpha_{j}}
\newcommand{\alfai}{\alpha_{i}}
\newcommand{\alfak}{\alpha_{k}}
\newcommand{\alfajpo}{\alpha_{j+1}}
\newcommand{\alfakpo}{\alpha_{k+1}}

\newcommand{\betzero}{\beta_{0}}
\newcommand{\betone}{\beta_{1}}
\newcommand{\bettwo}{\beta_{2}}
\newcommand{\betth}{\beta_{3}}
\newcommand{\betthr}{\beta_{3}}
\newcommand{\betfr}{\beta_{4}}
\newcommand{\betfv}{\beta_{5}}
\newcommand{\betsx}{\beta_{6}}
\newcommand{\betnmo}{\beta_{n-1}}
\newcommand{\betnmt}{\beta_{n-2}}
\newcommand{\betn}{\beta_{n}}
\newcommand{\betj}{\beta_{j}}
\newcommand{\betk}{\beta_{k}}
\newcommand{\betjpo}{\beta_{j+1}}
\newcommand{\betkpo}{\beta_{k+1}}

\newcommand{\gamzero}{\gamma_{0}}
\newcommand{\gamone}{\gamma_{1}}
\newcommand{\gamtwo}{\gamma_{2}}
\newcommand{\gamth}{\gamma_{3}}
\newcommand{\gamthr}{\gamma_{3}}
\newcommand{\gamfr}{\gamma_{4}}
\newcommand{\gamfv}{\gamma_{5}}
\newcommand{\gamsx}{\gamma_{6}}
\newcommand{\gamnmo}{\gamma_{n-1}}
\newcommand{\gamnmt}{\gamma_{n-2}}
\newcommand{\gamn}{\gamma_{n}}
\newcommand{\gamj}{\gamma_{j}}
\newcommand{\gamk}{\gamma_{k}}
\newcommand{\gamjpo}{\gamma_{j+1}}
\newcommand{\gamkpo}{\gamma_{k+1}}

\newcommand{\delzero}{\delta_{0}}
\newcommand{\delone}{\delta_{1}}
\newcommand{\deltwo}{\delta_{2}}
\newcommand{\delth}{\delta_{3}}
\newcommand{\delthr}{\delta_{3}}
\newcommand{\delfr}{\delta_{4}}
\newcommand{\delfv}{\delta_{5}}
\newcommand{\delsx}{\delta_{6}}
\newcommand{\delnmo}{\delta_{n-1}}
\newcommand{\delnmt}{\delta_{n-2}}
\newcommand{\deln}{\delta_{n}}
\newcommand{\delj}{\delta_{j}}
\newcommand{\delk}{\delta_{k}}
\newcommand{\deljpo}{\delta_{j+1}}
\newcommand{\delkpo}{\delta_{k+1}}

\newcommand{\xbar}{\bar{x}}
\newcommand{\ybar}{\bar{y}}
\newcommand{\dprime}{\prime\prime}
\newcommand{\rone}{r_{1}}
\newcommand{\rtwo}{r_{2}}

\large         
\begin{center}
ROOT FINDING TECHNIQUES THAT WORK
\vskip 0.5cm        
\normalsize
A. Melman \\
Department of Applied Mathematics \\ 
Santa Clara University  \\ 
Santa Clara, CA 95053 \\

\vskip 0.5cm                               
\vskip 0.5cm                               
\end{center}
\vskip 0.5cm                               

\begin{abstract}
In most introductory numerical analysis textbooks, the treatment of a single nonlinear equation often consists of a collection
of all-purpose methods that frequently do not work or are inefficient. These textbooks neglect to teach the importance
of adapting a method to the given problem,  and consequently also neglect to provide the tools to accomplish this.

Several general techniques are described here to incorporate the specific structure or properties of a nonlinear equation 
into a method for solving it. This can mean the construction of a method specifically tailored to the equation, or 
the transformation of the equation into an equivalent one for which an existing method is well-suited. The techniques
are illustrated with the help of several case studies taken from the literature. 
\end{abstract}
     
\vskip 1cm
\normalsize

%
%
\section{Introduction}
\label{introduction}

The numerical solution of a real nonlinear equation is frequently required in many areas of science and engineering, where it
often occurs as an important subproblem that needs to be solved repeatedly. There exist many well-known standard methods to achieve
this, typically iterative in nature, such as the well-known secant and Newton methods, and many others 
that can be found in any introductory numerical analysis textbook (Newton's method can even be found in high school calculus books). 
It is, in general, not easy to fully automate the solution of nonlinear
equations because of the many problems that may arise, such as, e.g., singularities or tightly clustered roots, to name but a few. The best situations
are those where the properties of the equation to be solved are well-known, as frequently occurs in specific applications. 
However, the inflexibility of standard methods often prevents efficient use of this information. As a general rule, even though it
may sometimes be easier said than done, every attempt should be made to tailor a method to the specific problem at hand, rather than using a 
general purpose method. The result will be a faster and more accurate method. An alternative approach is to transform the equation into
an equivalent one for which an already existing method is appropriate. Such considerations are unfortunately not commonly emphasized in textbooks. 

We will address these issues first by showing that several standard methods can be derived from a simple general principle that,
unlike these standard methods, has the 
flexibility to also generate methods that can take into account a particular problem's properties.
Secondly, we consider techniques to transform a given equation into an equivalent one with more useful properties, such as those
that guarantee the convergence of a particular method. Examples of the aforementioned principle and techniques can be found scattered 
throughout the literature, 
and we will use several of them as illustrative case studies. 

There are two main aspects to a numerical method: its construction and its subsequent convergence analysis. 
\emph{Global} convergence concerns global conditions under which 
the method is guaranteed to converge. Newton's method, for example, converges from any point to the right of the root of a convex increasing function, 
as its iterates are the roots of the tangents, which lie below such a convex function. Similar conditions can be obtained for other standard methods, 
although they are usually more complicated. \emph{Local} convergence of an iterative method is concerned with conditions that guarantee convergence
from a point that is sufficiently close to the root and also with the asymptotic rate at which the iterates converge. 
Denoting by~$x_{k}$ the iterates 
and by~$x^{*}$ the root to which they converge, a method is said to be of order~$q>0$ if
\bdis
\limsup_{k \rightarrow +\infty} \dfrac{|x_{k+1}-x^{*}|}{|x_{k}-x^{*}|^{q}} \leq C < +\infty \; . 
\edis
For~$q=1$, it is required that $C < 1$. The higher the value of~$q$, the faster the convergence. For example, when they converge to a simple root, 
i.e., when~$f^{\prime}(x^{*}) \neq 0$, then Newton's method (\cite[Ch.2]{Acton},\cite[Ch.3]{IK},\cite[Ch.5]{SB}) is of second order, the secant method 
(\cite[Ch.2]{Acton},\cite[Ch.3]{IK},\cite[Ch.5]{SB}) is of order $(1+\sqrt{5})/2$, and Halley's method (\cite{Salehov},\cite{ST}) is of third order.
Generally speaking, the more information (such as function and derivative values) is taken into account to construct the approximation the method 
is based on, the higher the order of convergence. An exhaustive treatment of methods and their convergence order can be found in~\cite{Traub}.

%
%
\section{Adaptive approximation and transformations}
\label{adaptapprox}

Consider the real nonlinear equation~$f(x)=0$ with $x \in \mathbb{R}$. To solve it, we formulate the following general adaptive approximation principle:
approximate~$f$ by another function~$g$ that satisfies certain requirements and for which the equation~$g(x)=0$ is easy to solve. 
An approximation to the root of~$f$ is then given by the (appropriate) root of the approximation~$g$.
An iterative method follows naturally from this principle by approximating~$f$ at a current iterate and generating the next iterate as a solution of $g(x)=0$.
There are two choices to make: the function~$g$ and the way in which it approximates~$f$. 
The approximation function is often (and
naturally) required to mimic the behavior of the function it approximates as much as possible. However, different requirements are sometimes imposed,
e.g., when the iterates need to be constrained in a certain way, as will be the case in Example~\ref{Pellet}. 

As an example, we now show how three well-known methods follow from this simple principle. To avoid interrupting the exposition, 
we will implicitly assume that
all expressions are valid, e.g, if a number appears in the denominator, it is assumed to be nonzero.
The first method chooses a line as an appropriate approximation function~$g$, 
i.e., $g(x)=\alpha + \beta x$, and the approximation conditions as requiring that $g$ coincide with~$f$ to first order, i.e., in function and first derivative 
values, at a given point~$\xbar$, so that 
\bdis
\left \{
\begin{array}{l} 
f(\xbar)=g(\xbar) = \alpha + \beta \xbar \\
f^{\prime}(\xbar) = g^{\prime}(\xbar) = \beta  \; ,
\end{array}
\right .
\edis
from which we obtain that $g(x) = f(\xbar) - \xbar f^{\prime}(\xbar) + f^{\prime}(\xbar) x $. The next iterate is the root of~$g$, given by 
$\xbar-f(\xbar) / f^{\prime}(\xbar)$, which is exactly one step of Newton's famous method that generates iterates $\{x_{k}\}$, $k \in \mathbb{N}$,
defined by the \emph{iteration formula}
\bdis
x_{k+1} = x_{k} - \dfrac{f(x_{k})}{f^{\prime}(x_{k})} \; . 
\edis
Of course, this is not a surprise because the Newton iterate
is often explained geometrically as the root of the tangent to~$f$, i.e., the root of the linear approximation to~$f$ at a certain point, 
which is precisely~$g$. 

The second method also picks~$g(x) = \alpha  + \beta x$, but changes the approximation conditions to the requirement that~$f$ and~$g$ coincide 
in function value at $\xbar$ and at one other (distinct) point $\ybar$, i.e., 
\bdis
\left \{
\begin{array}{l} 
f(\xbar) = g(\xbar) = \alpha + \beta \xbar  \\
f(\ybar) = g(\ybar) = \alpha + \beta \ybar \; ,
\end{array}
\right .
\edis
which means that the function~$g$ is found by linear interpolation at~$\bar{x}$ and~$\bar{y}$. This set of linear equations is easily solved for 
$\alpha$ and $\beta$, so that the next iterate, obtained from $-\alpha/\beta$, is given by
\bdis
\xbar - \dfrac{f(\xbar)(\xbar-\ybar)}{f(\xbar) - f(\ybar)} \; ,
\edis
which turns out to be one step of the secant method, namely, a method that can be viewed as obtained from Newton's method 
by approximating the derivative by a finite difference. Its iteration formula is given by
\bdis
x_{k+1} = x_{k} - \dfrac{f(x_{k})(x_{k}-y_{k})}{f(x_{k})-f(y_{k})} \; . 
\edis

The third method, Halley's method, the ``method of osculating hyperbolae'', can also be derived 
using the same principle we just illustrated (see, e.g., \cite{Salehov} and~\cite{ST}). 
In this case, $g(x) = \alpha + \beta/(x-\gamma)$, and the approximation requires that 
$f$ and $g$ coincide up to second derivatives at a given point, resulting in the iteration formula
\bdis
x_{k+1} = x_{k} - \dfrac{2f(x_{k})f^{\prime}(x_{k})}{2(f^{\prime}(x_{k}))^{2}-f(x_{k})f^{\prime\prime}(x_{k})} \; . 
\edis

Several more existing methods can be derived in this unifying way, but its main advantage lies in its flexibility. Consider, 
for example, the equation $f(x)=0$, where~$f$ is of the form $f(x)=f_{1}(x)+f_{2}(x)$. Standard methods are typically defined by an iteration 
formula rigidly applied to all of~$f$, based on an approximation function that does not necessarily bear any relation to~$f$. 
Instead, a more efficient method can often be obtained by approximating the functions $f_{1}$ and $f_{2}$ each in a different and more appropriate way.

However, it is not always possible or easy to find a convenient approximation. In such cases, one might attempt to adapt the problem to a certain
method, by which is meant formulating an equivalent problem for which a certain method exhibits desired convergence properties.
Such an equivalent problem might be obtained by a transformation of variables $x=\varphi(z)$ for a suitable function~$\varphi$ so that
the equivalent problem becomes $f(\varphi(z))=0$. Another possibility is the use of a nonzero multiplier $\mu(x)$ that transforms
the problem into $\mu(x)f(x)=0$, which has the same solutions as $f(x)=0$. Sometimes, a combination of techniques is called for.
In the following section, we take a detailed look at concrete applications of these techniques. 

%
%
\section{Case studies}
\label{cases}
We present four examples, taken from the literature, to illustrate the techniques of the previous section 
for both the construction and analysis of methods for solving nonlinear equations. In the first example, a method
is developed for solving an equation using adaptive approximation, where an approximation function is chosen to 
resemble the function that is being approximated. 

%
%
\begin{example}
\label{seceq_BNS}
{\bf (Secular Equation - adaptive approximation)}
This example considers the solution of a so-called \emph{secular equation} (\cite{GLB}), which lies at the heart of the fast and widely used divide and conquer method 
from~\cite{Cuppen} to compute the eigenvalues of a symmetric matrix.
After some simplification, this secular equation takes the form $f(x)=0$, where
\beq
\label{seceq}
f(x) \coloneqq 1 + \sum_{j=1}^{n} \dfrac{b_{j}}{d_{j} -x} \; ,
\eeq 
with $b_{j} > 0$ for all~$j$ and $d_{1} < d_{2} < \dots < d_{n-1} < d_{n}$. This function has $n$ poles and $n$ roots, 
one on each interval $(d_{j},d_{j+1})$ for $j=1,\dots,n-1$, and one on $(d_{n}, +\infty)$. 
The value of~$n$ can be very large and the goal is to compute all of the roots quickly and accurately. 

For the sake of this example, we concentrate on the $i$th root, with $1 \leq i \leq n-1$, after the origin is shifted to $d_{i}$, so that 
from here on we consider the computation of the unique root of~$f$ on the interval $(d_{i},d_{i+1}) = (0,d_{i+1})$, shown in Figure~\ref{figure_seceq_1}
as the lower curve in thicker line.
%
%
%
%
\begin{figure}[H]
\begin{center}
\raisebox{0ex}{\includegraphics[width=0.4\linewidth]{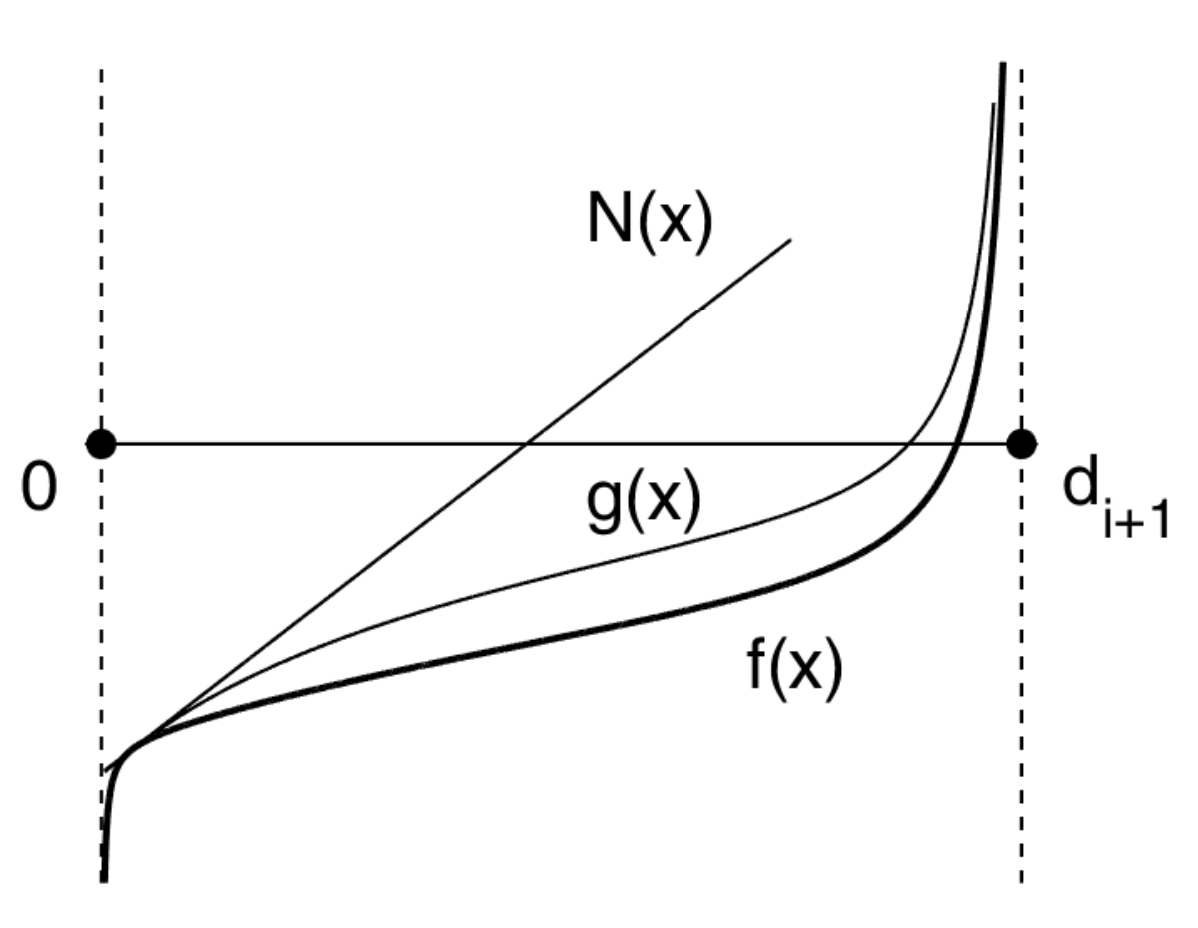}}
\caption{The functions~$f(x)$, $g(x)$, and~$N(x)$ on $(0,d_{i+1})$.}
\label{figure_seceq_1}            
\end{center}
\end{figure}
Figure~\ref{figure_seceq_1} illustrates a typical (although mild) situation, where a root's location is close 
to one of the singularities, as often happens in practice. A method such as Newton's method would generate iterates outside the interval 
from most starting points in $(0,d_{i+1})$, and trying to find a starting point to guarantee convergence might require a larger effort than the 
actual computation of the root. It is, in fact, not surprising that Newton's method would not be appropriate here as it is based on a linear approximation, 
while the function~$f$ is rational. We have made a comparison to Newton's method here simply because it is well-known, and similar reasoning applies to other 
standard methods. 
\newline {\bf \underline{Method construction.}} We can expect to obtain a more efficient method by using an approximation 
that more closely resembles the function~$f$ and distinguishes between different parts of~$f$. 
There are several candidates for such an approximation, which were extensively studied in~\cite{BNS},
\cite{RCL}, \cite{Melman_SecEq1}, and \cite{Melman_SecEq2}. Here we will consider the approximation used in~\cite{BNS}.
To do this, we first rewrite the function~$f$ 
as $f=1+f_{1}+f_{2}$, with    
\beq
\nonumber
f_{1}(x) \coloneqq \sum_{j=1}^{i} \dfrac{b_{j}}{d_{j} -x} \;\;\;\; \text{and} \;\;\;\; f_{2}(x) \coloneqq \sum_{j=i+1}^{n} \dfrac{b_{j}}{d_{j} -x} \; , 
\eeq
and construct an iterative numerical method by separately approximating~$f_{1}$ and~$f_{2}$ by~$g_{1}$ and~$g_{2}$, respectively, so that, at a given point, 
$g_{1}$ agrees with $f_{1}$ in function and first derivative values, while $g_{2}$ does the same for $f_{2}$.
The approximation~$g$ to~$f$ is then defined as $g \coloneqq 1+g_{1}+g_{2}$.
One possible choice, as in~\cite{BNS}, is to use the rational functions  
\beq
\nonumber
g_{1}(x) \coloneqq \dfrac{\alpha}{\beta -x} \;\;\;\; \text{and} \;\;\;\; g_{2}(x) \coloneqq \gamma + \dfrac{\delta}{d_{i+1} - x} \; .  
\eeq
In other words, at a given point~$\xbar$, the following conditions must be satisfied:
\beq
\label{approxcond}
\left \{ \begin{array}{l}
         f_{1}(\xbar)=g_{1}(\xbar) \\
         f_{1}^{\prime}(\xbar)=g_{1}^{\prime}(\xbar) \\
         \end{array}
\right.
\;\;\;\;\;\; \text{and} \;\;\;\;\;\;
\left \{ \begin{array}{l}
         f_{2}(\xbar)=g_{2}(z) \\
         f_{2}^{\prime}(\xbar) = g_{2}^{\prime}(\xbar) \; . \\
         \end{array}
\right.
\eeq 
The conditions in~(\ref{approxcond}), which require the same computational effort as Newton's method (function and derivative values),
determine the parameters $\alpha$, $\beta$, $\gamma$, and $\delta$, which, in turn, define the functions $g_{1}$ and~$g_{2}$.
The approximation~$g$ is well defined on the interval $(0,d_{i+1})$ since the approximation conditions yield
\bdis
\beta = \xbar + \dfrac{f_{1}(\xbar)}{f_{1}^{\prime}(\xbar)} 
= \dfrac{\xbar f_{1}^{\prime}(\xbar)+f_{1}(\xbar)}{f_{1}^{\prime}(\xbar)} 
= \dfrac{1}{f_{1}^{\prime}(\xbar)} \sum_{j=1}^{i} \dfrac{b_{j}d_{j}}{(d_{j} - \xbar)^2}  < 0 \; .
\edis
Moreover, $\alpha=(\beta-\xbar)^{2}f_{1}^{\prime}(\xbar) > 0$ and $\delta=(d_{i+1}-\xbar)^{2}f_{2}^{\prime}(\xbar) > 0$, so that $g^{\prime}(x) >0$
on the interval, implying that~$g$ is strictly increasing (to~$+\infty$). It must therefore have a unique root on~$(0,d_{i+1})$ if~$\xbar$ is chosen 
such that $f(\xbar)=g(\xbar)<0$.
The next iterate of a method based on this approximation is the root of~$g$ in the interval $(0,d_{i+1})$, which is obtained 
as the appropriate solution of a simple quadratic, since
\bdis
1 + \dfrac{\alpha}{\beta -x} + \gamma + \dfrac{\delta}{d_{i+1} - x} = 0 
\Longrightarrow 
(1+\gamma)(\beta-x)(d_{i+1}-x)+\alpha(d_{i+1}-x)+\delta(\beta-x) = 0 \; .  
\edis
The approximation at a point~$\xbar$, deliberately chosen to be far from the root, is shown in Figure~\ref{figure_seceq_1} in a thinner line just above the 
function~$f$. For comparison, we have also shown the tangent~$N(x)$ at the same point~$\xbar$, which clearly demonstrates the advantage of 
a rational approximation. 

\noindent {\bf \underline{Convergence properties.}} To address the global convergence of the method just described, 
we now first show that~$g$ dominates~$f$ on~$(0,d_{i+1})$. 
Since $g_{1}(x) = \alpha/(\beta-x)$ approximates~$f_{1}(x)$ to first order at~$\xbar \in (0,d_{i+1})$, 
$(\beta -x)/\alpha$ approximates $1/f_{1}(x)$ to first order at $\xbar$, i.e., $(\beta-x)/\alpha$ is the linear approximation to $1/f_{1}$
at~$\xbar$. A straightforward calculation yields
\beq
\label{ddfoneinv}
\lb \dfrac{1}{f_{1}}\rb^{\dprime} =\dfrac{2(f_{1}^{\prime})^2-f_{1}f_{1}^{\dprime}}{f_{1}^{3}} 
=  \dfrac{-2(-f_{1}^{\prime})^2+(-f_{1})(-f_{1}^{\dprime})}{(-f_{1})^{3}} \; . 
\eeq
The function~$-f_{1}$ is positive on the interval~$(0,d_{i+1})$ and it satisfies the conditions of Lemma~2.3 in~\cite{Melman_SecEq1} for $\rho=-1$, 
a parameter in that lemma, which in this case states that $-2(-f_{1}^{\prime})^2+(-f_{1})(-f_{1}^{\dprime}) \geq 0$. It then follows from~(\ref{ddfoneinv}) 
that $1/f_{1}$ is a convex function on~$(0,d_{i+1})$, so that it dominates its linear approximation at any point in that interval. As a result,          
\bdis
\dfrac{\beta - x}{\alpha} \leq \dfrac{1}{f_{1}(x)} \Longrightarrow
\dfrac{\alpha}{\beta - x} \geq f_{1}(x) \; ,
\edis
i.e., $g_{1}(x) \geq f_{1}(x)$. We also have that $g_{2}(x) = \gamma +  \delta/(d_{i+1}-x)$ approximates~$f_{2}$ to first order at~$\xbar$,
which means that $\gamma(d_{i+1}-x) + \delta$ is the linear approximaton of $(d_{i+1}-x)f_{2}(x)$. Some algebra yields
\bdis
(d_{i+1}-x)f_{2}(x) = \lb \sum_{j=i+1}^{n} b_{j} \rb - \sum_{j=i+1}^{n} \dfrac{b_{j}(d_{j}-d_{i+1})}{d_{j}-x} \; , 
\edis 
which is a concave function because $b_{j}(d_{j}-d_{i+1}) \geq 0$ when $j \geq i+1$. 
This means that it is dominated by its linear approximation, leading to 
\bdis
\gamma(d_{i+1}-x) + \delta \geq (d_{i+1}-x)f_{2}(x) \Longrightarrow g_{2}(x) = \gamma +  \dfrac{\delta}{d_{i+1}-x} \geq f_{2}(x) \; .
\edis
As a result, we obtain that $g(x) \geq f(x)$ on~$(0,d_{i+1})$. Consequently, if~$\xbar$ lies to the left of the root, then $f(\xbar)=g(\xbar) <0$,
implying that the unique root of~$g$, which necessarily lies to the right of~$\xbar$, also lies to the left of the root of~$f$, ensuring monotonic 
convergence. Moreover, it was shown in ~\cite{Melman_SecEq1} that the convergence order is quadratic.

We conclude by deriving a starting point~$x_{0} \in (0,d_{i+1})$ satisfying $f(x_{0}) < 0$. Such a point can be found by observing that 
\bdis
1 + \sum_{j=1}^{i-1} \dfrac{b_{j}}{d_{j} - d_{i+1}} - \dfrac{b_{i}}{x} + \dfrac{b_{i+1}}{d_{i+1}-x} + \sum_{j=i+2}^{n} \dfrac{b_{j}}{d_{j} - d_{i+1}}  
\geq f(x) \; ,
\edis
so that initial point~$x_{0}$ can be found as the root in~$(0,d_{i+1})$ of the strictly increasing function
\bdis
\lb 1 + \sum_{\stackrel{j=1}{j \neq i,i+1}}^{n} \dfrac{b_{j}}{d_{j} - d_{i+1}} \rb - \dfrac{b_{i}}{x} + \dfrac{b_{i+1}}{d_{i+1}-x} \; ,
\edis
which is obtained as the appropriate root of a quadratic.

An additional consideration in the case of a parallel computation of the~$n$ roots of~$f$ is that the efficiency of such an implementation 
is determined by the root that requires the most time, so that uniform performance of the method is also important. In practice, the more an 
approximation resembles the function, the better this requirement will be fulfilled. 
\end{example}

In the following example, the equation from Example~\ref{seceq_BNS} is solved in a different way, namely, by first carrying out a transformation
of variables before applying adaptive approximation with a function having similar properties as the one being approximated.

%
%
\begin{example}
\label{seceq_Bakchik}
{\bf (Secular Equation - transformation and adaptive approximation)}
In Example~\ref{seceq_BNS} a method was derived to compute the root of~$f$, defined in~(\ref{seceq}), on the interval $(0,d_{i+1})$ by using rational 
approximations. 
Here, the approach from~\cite[Section 3.3]{Melman_SecEq1} is used, which consists of transforming the variable to obtain an equivalent, 
but more convenient, equation.

\noindent {\bf \underline{Method Construction.}} 
A relatively natural idea is to try and mitigate the effect of the singularities at the endpoints of the interval. 
A transformation that then suggests itself (\cite{Melman_SecEq1}), is to set $x=1/z$, which eliminates the singularity at the origin by sending it 
to infinity, so that the original interval is mapped to $(1/d_{i+1},+\infty)$, and the problem becomes the computation of the (necessarily) unique root of
$f(1/z)=0$ on this interval. For convenience, we define $F(z)=f(1/z)$.
So far, this is an idea that seems reasonable, but we still need to show that the resulting transformed equation $F(z)=0$
exhibits properties that make it easier to solve than $f(x)=0$. Straightforward algebra shows that 
\beq
\nonumber           
F(z) = 1 + \sum_{j=1,j \neq i}^{n} \dfrac{b_{j}}{d_{j}}-b_{i}z +  \sum_{j=1,j \neq i}^{n} \dfrac{b_{j}/d_{j}^{2}}{z-1/d_{j}}  \; .
\eeq                
All the negative singularities of $F$ lie in the interval $(1/d_{i-1},0)$ and all its positive singularites lie in the interval $(0,1/d_{i+1})$. 
Its first and second derivatives show that that $F$ is strictly decreasing and convex on the interval $(1/d_{i+1},+\infty)$, with 
\bdis
\lim_{z \rightarrow (1/d_{i+1})^{+}} F(z) = +\infty \qquad \text{and} \qquad \lim_{z \rightarrow +\infty} F(z) = -\infty \; .
\edis
Because~$F$ is convex, it dominates its linear approximation, so that Newton's method converges monotonically from any starting point 
between $1/d_{i+1}$ and the root. Figure~\ref{figure_seceq_2} shows the function~$F(z)$ (top curve in thicker line), along with its tangent~$N(z)$ 
at a particular point. The root of the tangent is the next Newton iterate from that point.

The transformation of variables $x=1/z$ has transformed the original equation into one for which Newton's method is guaranteed to converge, 
which was not the case for the original equation, as we saw in Example~\ref{seceq_BNS}. However, the equation still involves a rational function that 
should preferably be approximated by a rational, rather than a linear, function. To construct such an approximation, we observe that the behavior 
of~$F$ on $(1/d_{i+1},+\infty)$ is significantly affected by its singularity at $1/d_{i+1}$, which indicates that, if possible, the singularity should 
be included in the approximation. Writing~$F$ as   
\beq
\nonumber
F(z) = F_{1}(z) + \dfrac{b_{i+1}/d_{i+1}^{2}}{z-1/d_{i+1}} \; ,
\eeq
we approximate it at a point $\xbar$ by
\beq
\nonumber
G(z) = \alpha + \beta z + \dfrac{b_{i+1}/d_{i+1}^{2}}{z-1/d_{i+1}} \; ,
\eeq
where $\alpha + \beta z$ is the linear approximation of~$F_{1}$ at $\xbar$. Since~$\beta=F_{1}^{\prime}(\xbar) < 0$ for any~$\xbar \in (1/d_{i+1},+\infty)$, 
$G$ is strictly decreasing, while it becomes unbounded as  $z \rightarrow (1/d_{i+1})^{+}$, implying that it has a unique root on $(1/d_{i+1},+\infty)$. 
The next iterate of
a method based on this approximation is therefore the root of~$G$, which is found as the appropriate root of the quadratic 
$(z-1/d_{i+1})(\alpha + \beta z)+b_{i+1}/d_{i+1}^{2}$.
%
%
%
%
\begin{figure}[H]
\begin{center}
\raisebox{0ex}{\includegraphics[width=0.4\linewidth]{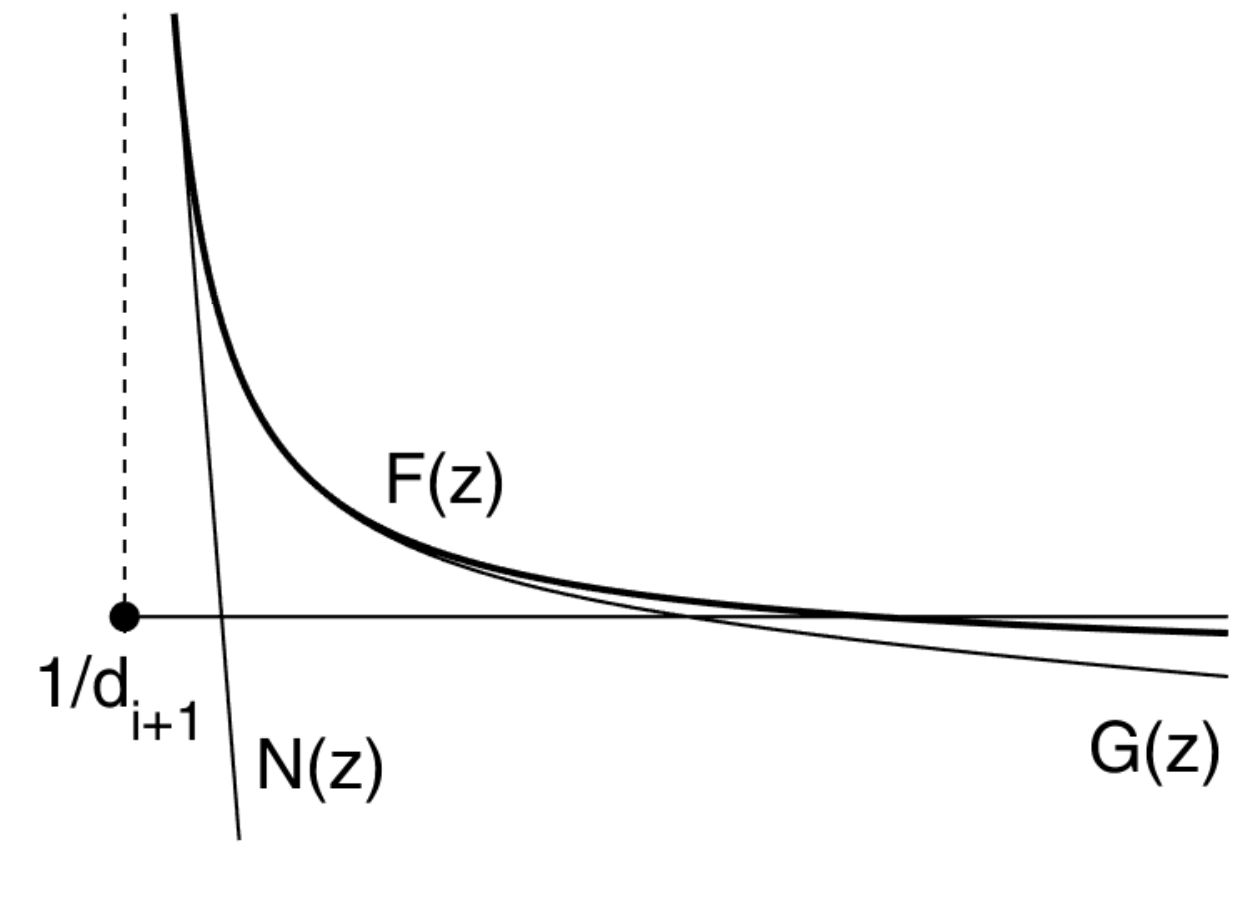}}
\caption{The functions~$F(z)$, $G(z)$, and~$N(z)$ on $(1/d_{i+1},+\infty)$.}
\label{figure_seceq_2}            
\end{center}
\end{figure}

\noindent {\bf \underline{Convergence properties.}} Since, like~$F$, $F_{1}$ is convex, $\alpha + \beta z$ will be dominated by~$F_{1}$, 
and the entire approximation will therefore be dominated by~$F$. As a result, the convergence of a method based on this approximation will be, like Newton's 
method, monotonic from any starting point between $1/d_{i+1}$ and the root of~$F$. Unlike Newton's method, it is globally convergent from any point in 
$(1/d_{i+1},+\infty)$: if the initial iterate lies to the right of the root of~$F$, the second iterate will lie between $1/d_{i+1}$ and the root.
Moreover, because the approximation resembles~$F$ more closely than a line, we expect such a method to perform better than Newton's method,
while requiring the same computational effort, namely, the computation of~$F(\bar{z})$ and~$F^{\prime}(\bar{z})$. 
Figure~\ref{figure_seceq_2} clearly shows this to be the case, where the nonlinear approximation (at the same point
that was used for the linear approximation) is the curve just below~$F(z)$, drawn in thinner line. Its root is a much better approximation to the root of~$F$
than the one obtained from Newton's method, which is the root of the linear approximation (tangent line). The quadratic convergence of such a method
immediately follows from the quadratic convergence of Newton's method. 

Similar arguments as in Example~\ref{seceq_BNS} show that a initial point to the right of the root of~$F$ is obtained as the root of
\bdis
1 + \sum_{j=1,j \neq i}^{n} \dfrac{b_{j}}{d_{j}} +  \sum_{j=1,j \neq i,i+1}^{n} \dfrac{b_{j}/d_{j}^{2}}{1/d_{i+1}-1/d_{j}}  
- b_{i}z + \dfrac{b_{i+1}/d_{i+1}^{2}}{z-1/d_{i+1}} \; , 
\edis
which dominates~$F$, while a point to the left of the root is obtained as the root of
\bdis
1 + \sum_{j=1,j \neq i}^{n} \dfrac{b_{j}}{d_{j}} - b_{i}z + \dfrac{b_{i+1}/d_{i+1}^{2}}{z-1/d_{i+1}} \; , 
\edis
which is dominated by~$F$. Both are computed as the appropriate root of a quadratic. An initial point can then be chosen as the average of these two points.
In practice, the performance of the method in this example is similar to that of the method in Example~\ref{seceq_BNS}.

We conclude with the observation that the method we obtained can easily be improved by including more terms of~$F$ in the approximation function~$G$.
For example, we could define 
\beq
\nonumber
G(z) = \alpha + \beta z + \dfrac{b_{i+1}/d_{i+1}^{2}}{z-1/d_{i+1}} + \dfrac{b_{i+2}/d_{i+2}^{2}}{z-1/d_{i+2}} 
                   + \dfrac{b_{i+3}/d_{i+3}^{2}}{z-1/d_{i+3}} + \dfrac{b_{i+4}/d_{i+4}^{2}}{z-1/d_{i+4}} \; ,
\eeq
and compute its root with a method like the one just obtained. Such an approach could be advantageous for large values of~$n$, as is often the case in practice.
\end{example}

In the next example a nonzero multiplier is used to transform the equation into one for which Newton's method exhibits global convergence.

%
%
\begin{example}
\label{knapsack}   
{\bf (Knapsack - multiplier)}
The equation we consider in this example is obtained in the course of solving the \emph{nonlinear continuous knapsack problem} from~\cite{Melman_NLK}.
In the original and simplest form of the (discrete and linear) knapsack problem, a knapsack of limited volume is filled with items that have different 
volumes and values, 
with the goal of maximizing the total value of the items in the knapsack. It has been generalized to nonlinear continuous problems with many applications 
bearing no relation to knapsacks (see~\cite{H}), as in our case here, where the problem originated in the scheduling of the servicing of chemical processing units.

Before we state the equation to be solved, we define $h(x)= 1-(1+1/x)e^{-1/x}$ on~$(0,+\infty)$ and its inverse function $\varphi(x)=h^{-1}(x)$, 
defined on $(0,1)$. We note that~$\varphi$ is well-defined since $h'(x) = -x^{-3}e^{-1/x} < 0$ for $x>0$. Figure~\ref{figure_hhinv} shows the 
functions~$h$ and $\varphi$. The function~$\varphi$ will play a crucial rule, and to gain a better understanding of its properties, we make the following 
observations. Set $y=\varphi(x)$ to obtain                 
\bdis
y=\varphi(x) \; \Longrightarrow \; h(y)=x \; \Longrightarrow \; h^{\prime} (y) y^{\prime} = 1 \; \Longrightarrow \; y^{\prime} = -y^{3}e^{1/y} \; ,
\edis
from which it follows that $y^{\prime\prime} = -y(3y-1)e^{1/y} y^{\prime} = y^{4}(3y-1)e^{2/y}$. These calculations show that~$\varphi$ is
a strictly decreasing function with a single inflection point at $h(1/3)$. Moreover, 
\bdis
\lim_{x \rightarrow 0^{+}} \varphi(x) = +\infty
\;\; \text{and} \;\;
\lim_{x \rightarrow 1^{-}} \varphi^{\prime}(x) = -\infty \; .
\edis
The function~$\varphi$ does not have an explicit functional expression: to compute~$y=\varphi(x)$, one needs to compute the solution~$y$ of the nonlinear 
equation $h(y)=x$, which will briefly be addressed further on. 
%
%
%
%
\begin{figure}[H]
\begin{center}
\raisebox{0ex}{\includegraphics[width=0.35\linewidth]{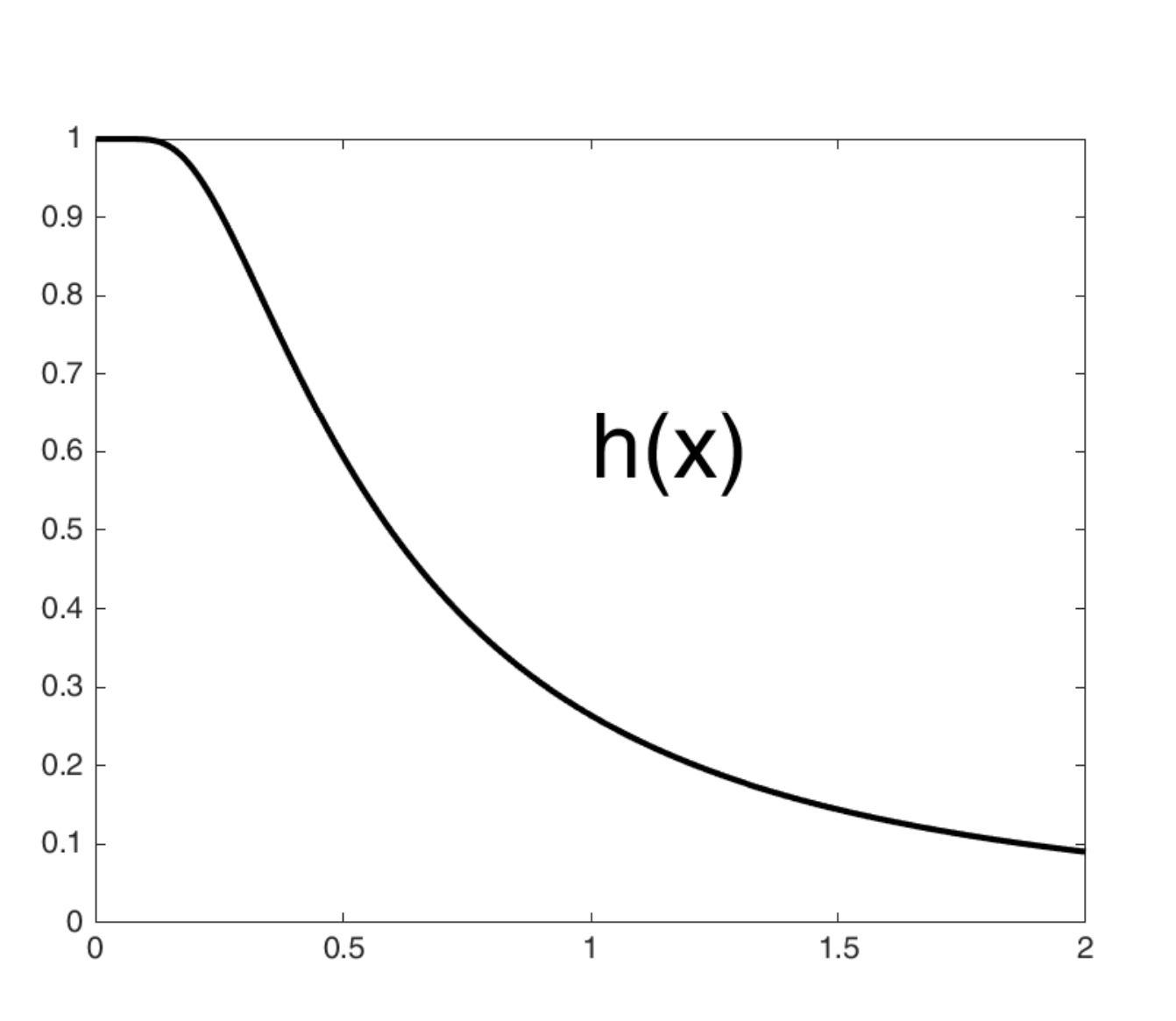}}
\hskip 2.0cm
\raisebox{0ex}{\includegraphics[width=0.35\linewidth]{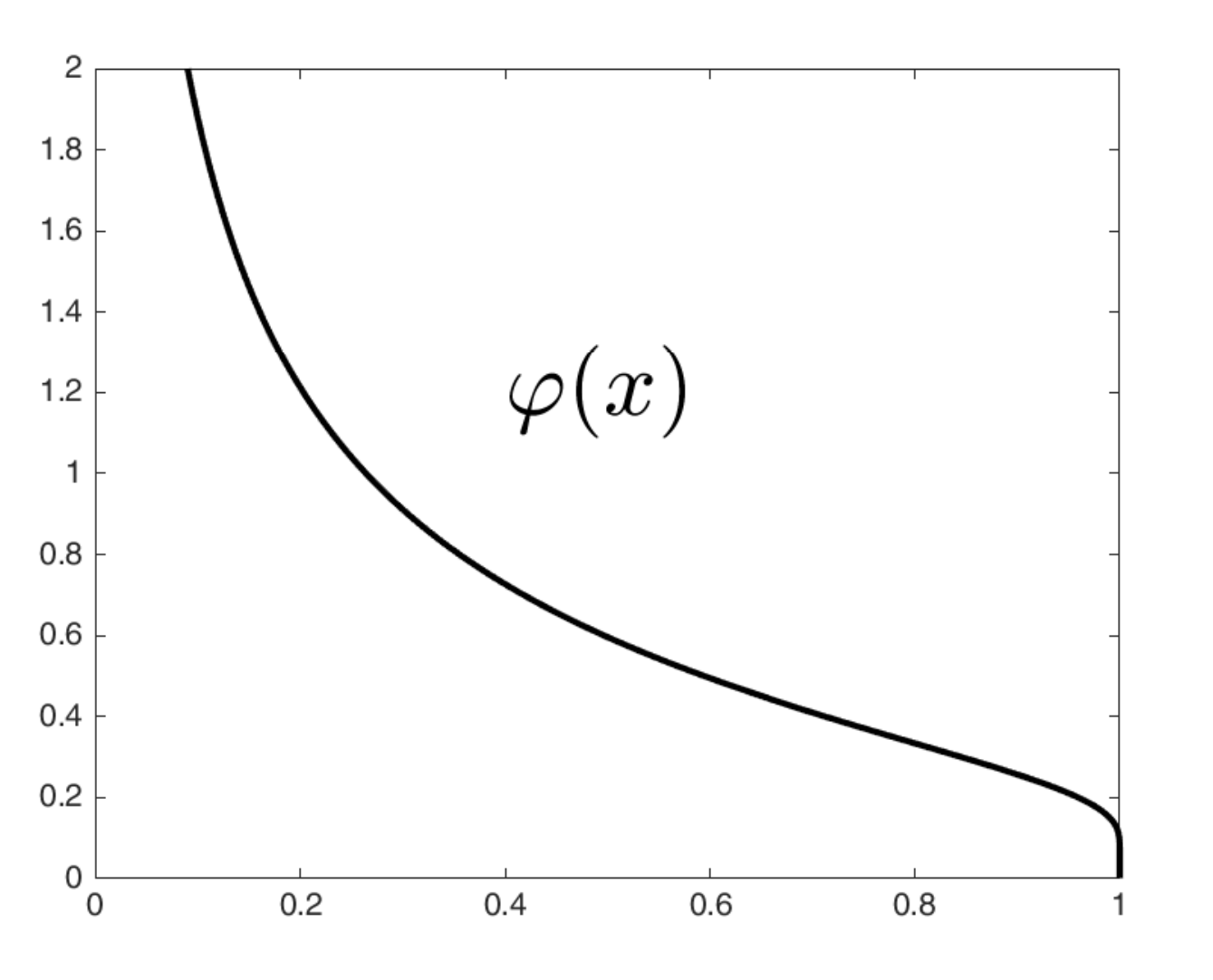}}
\caption{The functions~$h(x)$ and~$\varphi(x)$.}
\label{figure_hhinv}            
\end{center}
\end{figure}
We are now ready to state the problem to be solved, which is to compute the root of the function~$f$, originating from the dual formulation of the problem,
defined by
\bdis
f(x) \coloneqq \sum_{j=1}^{n} \alfaj \varphi(\betj x) - K \; , 
\edis
where $\alfaj,\betj,K>0$, on the interval $(0,\gamma)$, with $\gamma \coloneqq \min_{j} \{ 1/\betj \}$. Because it is a positive linear combination of scaled
versions of~$\varphi$, the properties of~$f$ are similar to those of~$\varphi$~: it is a strictly decreasing function, which becomes unbounded as 
$x \rightarrow 0^{+}$, and has an unbounded derivative when $x \rightarrow \gamma^{-}$. It is the curve in thicker line in both graphs of Figure~\ref{figure_ffion}. 
In what follows, we assume that $K > \lim_{x \rightarrow \gamma^{-}} f(x)$, implying that~$f$ has a unique root on $(0,\gamma)$.
%
%
%
%
\begin{figure}[H]
\begin{center}
\raisebox{0ex}{\includegraphics[width=0.35\linewidth]{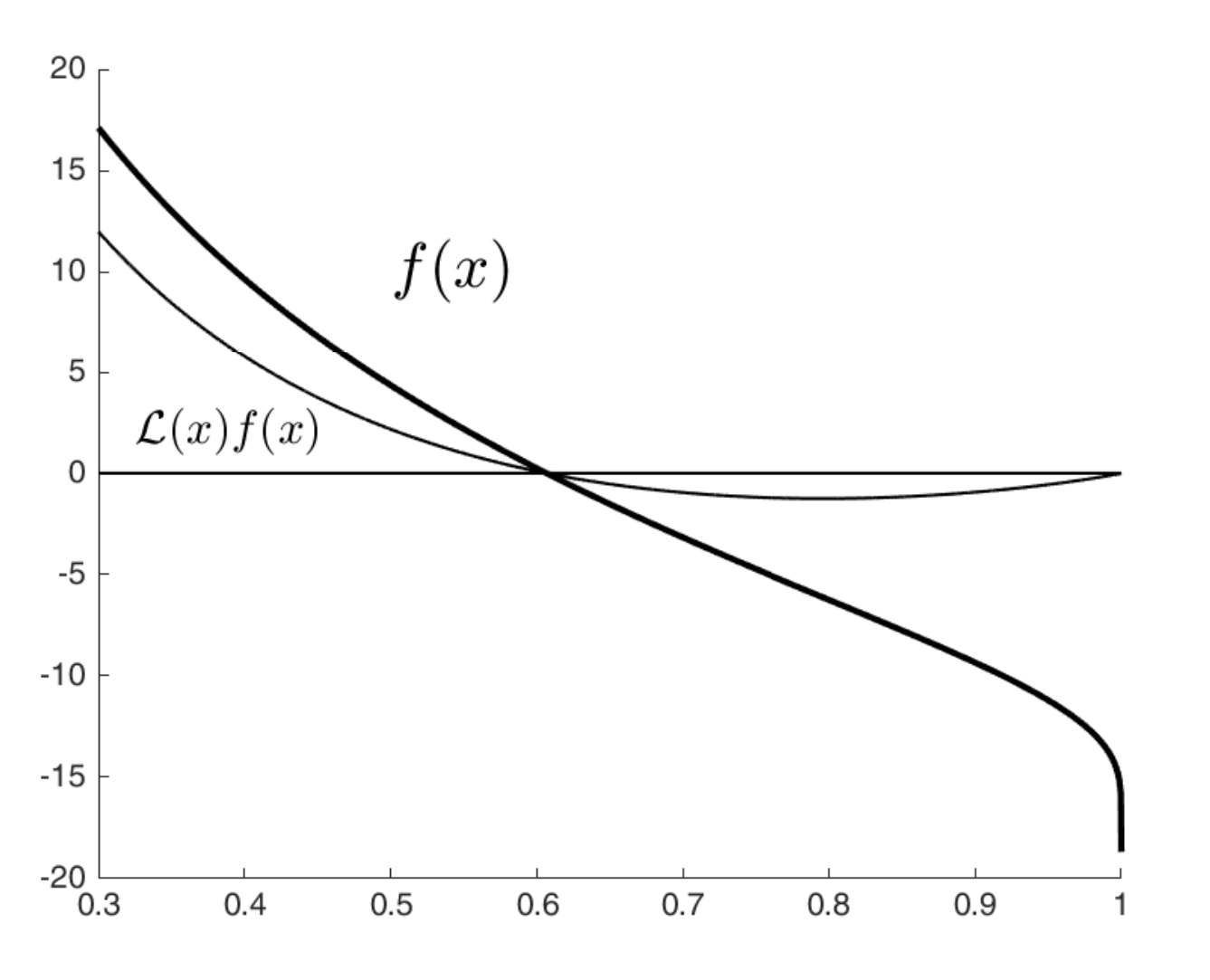}}
\hskip 1.0cm
\raisebox{0ex}{\includegraphics[width=0.35\linewidth]{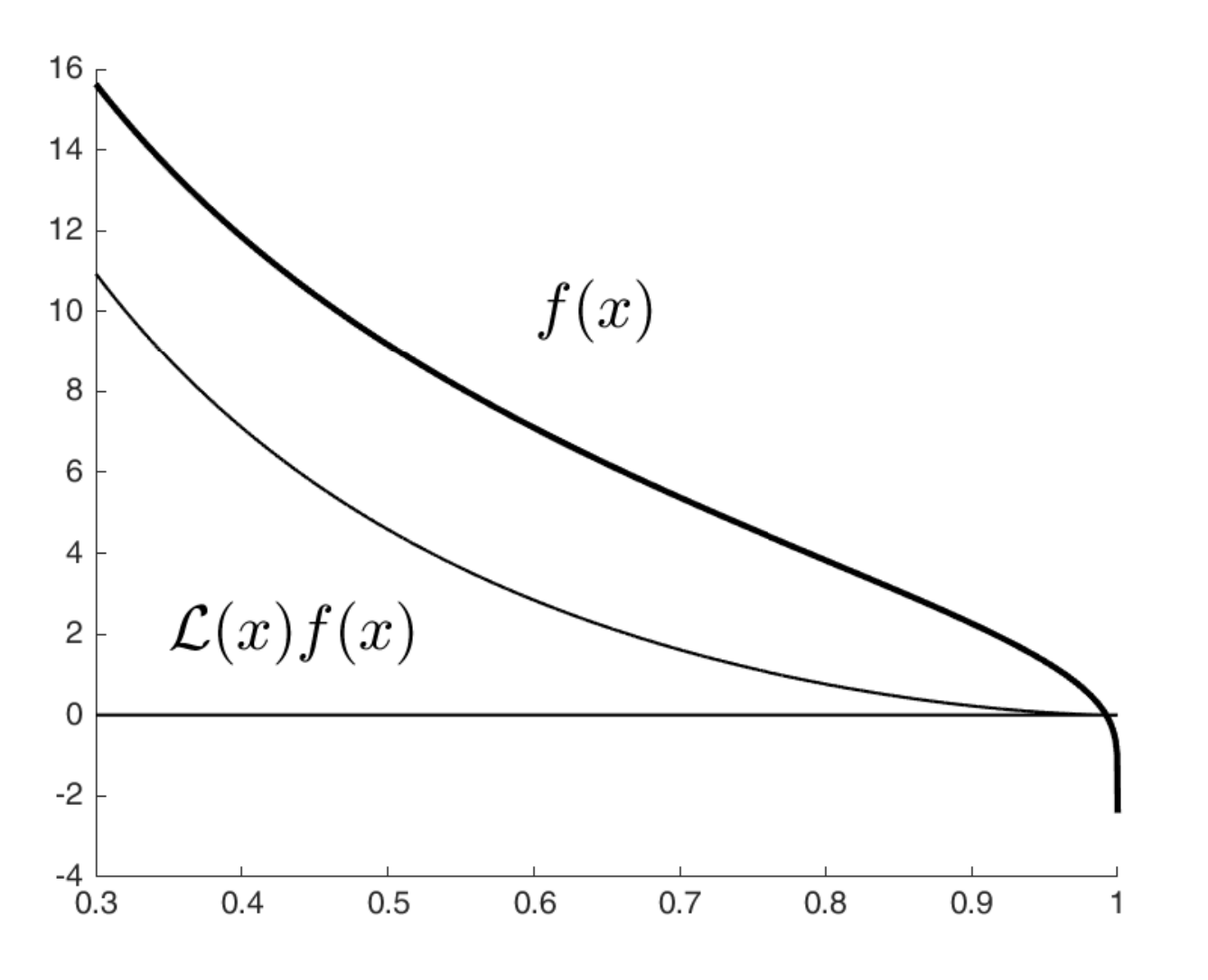}}
\caption{The functions~$f(x)$ and~$\mathcal{L}(x)f(x)$.}
\label{figure_ffion}            
\end{center}
\end{figure}
\noindent {\bf \underline{Convexifying multiplier.}} 
The function~$f$ exhibits several difficulties: it does not have an explicit functional expression, and it does not have a simple shape, as it is neither 
convex nor concave over~$(0,\gamma)$. The graph on the left in Figure~\ref{figure_ffion} shows a mild situation, whereas, on the right, the root lies very 
close to~$\gamma$, as inevitably happens in practice. In the latter case, a standard method (including Newton's method) can easily produce iterates that fall outside the interval, 
unless a starting point is available that already lies very close to the root. There does not seem to be an obvious appropriate approximation function that one 
could use, nor is a transformation of variables of much help. The best course of action may therefore be to try and use a multiplier to obtain a function with 
the same root, but with properties that make it easy to apply a standard method. In other words, we will adapt the problem to a method, rather than constructing
a specific method for the problem. An easy choice for such a standard method would be Newton's method, which is simple, asymptotically fast,
and guarantees convergence for convex and concave functions. 
The problem then becomes to find a nonzero multiplier~$\mathcal{L}(x)$ such that the function $\mathcal{L}(x)f(x)$ is either convex or concave.
In~\cite{Melman_NLK} it is shown that~$\mathcal{L}f$ is convex on~$(0,\gamma)$ if~$\mathcal{L}$ satisfies a specific differential inequality. 
The simplest such multiplier is $\mathcal{L}(x)=1-\gamma^{-1}x$. We will not reproduce its derivation as it is rather technical and irrelevant for our purposes,
but we will verify that $(1-\gamma^{-1}x)f(x)$ is indeed convex on~$(0,\gamma)$. To do this, we set $y=\varphi(\betj x)$, where $1 \leq j \leq n$, so that
$h(y)=1-(1+1/y)e^{-1/y}=\betj x$, from which it follows that
\beq
\label{expy}
e^{1/y} = \dfrac{y+1}{y(1-\betj x)} \; .
\eeq
We now use the expressions for $\varphi^{\prime}$ and $\varphi^{\prime\prime}$ that were previously derived, the expression for~$e^{1/y}$ from~(\ref{expy}),
and the fact that $\gamma^{-1}=\max_{j} \{ \betj \}$, to obtain that
\begin{eqnarray}
\lb (1- \gamma^{-1}x)\varphi(\betj x) \rb ^{\prime \prime} 
& = & -2\gamma^{-1}\betj \varphi^{\prime}(\betj x) + \betj^{2}(1-\gamma^{-1}x)\varphi^{\prime\prime}(\betj x) \nonumber \\
& = & 2\gamma^{-1}\betj y^{3}e^{1/y} + \betj^{2}(1-\gamma^{-1}x) y^{4}(3y-1) e^{2/y}   \nonumber  \\
& = & \betj y^{3}e^{1/y} \lb 2 \gamma^{-1} + \betj(1-\gamma^{-1}x) y(3y-1) e^{1/y} \rb  \nonumber  \\
& = & \betj y^{3}e^{1/y} \lb 2 \gamma^{-1} + \dfrac{\betj(1-\gamma^{-1}x)}{1-\betj x} (3y-1)(y+1) \rb  \nonumber  \\
& = & \betj y^{3}e^{1/y} \lb \dfrac{\betj(1-\gamma^{-1}x)}{1-\betj x} (3y^{2}+2y) 
                                                 +  2\gamma^{-1} - \lb \dfrac{1- \gamma^{-1}x}{1-\betj x} \rb \betj \rb \nonumber  \\
& \geq & \betj y^{3}e^{1/y} \lb \dfrac{\betj(1-\gamma^{-1}x)}{1-\betj x} (3y^{2}+2y) +  \gamma^{-1}  \rb   > 0 \; . \label{cnvxty}
\end{eqnarray}
Since
\bdis
(1-\gamma^{-1}x)f(x) = \sum_{j=1}^{n} \alfaj (1-\gamma^{-1}x)\varphi(\betj x) - K(1-\gamma^{-1}x)  
\edis
and $(1-\gamma^{-1}x)^{\prime\prime}=0$, we conclude with the help of~(\ref{cnvxty}) that~$\mathcal{L}f$ is convex on~$(0,\gamma)$.
Figure~\ref{figure_ffion} shows the function~$(1-\gamma^{-1}x)f(x)$ in thinner line. 
\newline {\bf \underline{Convergence properties.}} The convexity of~$\mathcal{L}f$ implies that if Newton's method 
is started from any point to the left of its root, then the iterates will converge monotonically to that root. Moreover, it is not hard to
show that the root is simple, so that Newton's method converges quadratically. 

A starting point to the left of the root of~$f$ can be found using the fact that, for $x>0$, $\gamma^{-1} x \geq \betj x$, so that
$\varphi(\gamma^{-1} x) \leq \varphi(\betj x)$, which in turn implies that 
\beq
\label{stptineq}
\sum_{j=1}^{n} \alfaj \varphi(\betj x) - K \geq \sum_{j=1}^{n} \alfaj \varphi(\gamma^{-1} x) - K \; .
\eeq    
The function in the right-hand side of~(\ref{stptineq}) is strictly decreasing, becomes unbounded at the origin, and is negative as $x \rightarrow 1^{-}$.
It must therefore have a unique root in~$(0,1)$ and, since~(\ref{stptineq}) shows that it is dominated by~$f$, that root must lie to the left of the root of~$f$.
A starting point $x_{0}$ is therefore given by
\bdis
x_{0} = \gamma \, h \hskip -0.1cm \lb \dfrac{K}{\sum_{j=1}^{n} \alfaj} \rb \; .
\edis

Although it is not the focus of this example, we conclude by briefly mentioning the function value computation of~$\varphi$.
From~$y=\varphi(x)$, we obtain $h(y)=1-1(1+1/y)e^{-1/y}= x$, so that~$y$ is the solution of a nonlinear equation on $(0,1)$. 
Setting $y=1/z$ transforms this equation into $1-(1+z)e^{-z}-x=0$, for which it was shown in~\cite{Melman_NLK} that Halley's 
method converges from any point in~$(0,1)$.
\end{example}

In the following example, a transformation of variables is used to facilitate adaptive approximation. In this case, 
the approximation function is chosen to satisfy constraints on the iterates, rather than to resemble the function it approximates.

%
%
\begin{example}
\label{Pellet}      
{\bf (Pellet - transformation and adaptive approximation)}
In this example, we consider Pellet's theorem (\cite[Th. (2,8)]{Marden}) for a polynomial $p(z)=\sum_{j=0}^{n} a_{j}z^{j}$. 
It states that if, for some $\ell$ with $1 \leq \ell \leq n-1$, $a_{\ell} \neq 0$, the real polynomial
$
q(z) \coloneqq |\an| z^{n} + \dots + |a_{\ell+1}| z^{\ell+1} - |a_{\ell}| z^{\ell} + |a_{\ell-1}|z^{\ell-1} + \dots  + |a_{0}| 
$
has two distinct positive roots, namely, the \emph{Pellet $\ell$-radii} $\rho_{1}$ and $\rho_{2}$, with $\rho_{1} < \rho_{2}$, 
then $p$ has exactly $\ell$ zeros in the closed disk $|z| \leq \rho_{1}$, and no zeros in the open annulus $\rho_{1} < |z| < \rho_{2}$.
In other words, the theorem can sometimes detect gaps between the moduli of zeros. It is a direct consequence of Rouch\'{e}'s theorem (\cite[Theorem 1.6]{Lang}).
The theorem has been generalized to matrix polynomials (\cite{BiniNoferiniSharify},\cite{Melman_GenPellet}), 
leading to a real polynomial of the same kind as~$q$. The graph of~$q$         
has a form that is very similar to the lower curve on the left in Figure~\ref{figure_trinom} in thicker line. It is of overriding importance 
for the equation $q(x)=0$ to be solved by a method whose iterates converge to the solutions from the inside of the interval $[\rho_{1},\rho_{2}]$, 
because, in such a case, one has the option to stop iterating at any moment and still have correct bounds. 
If, on the other hand, the iterates converge from outside the interval, 
none of the iterates provide correct bounds until they have fully converged.

An efficient method to compute the Pellet radii from inside the interval was derived in~\cite{Melman_Pellet}, and consists of two phases. 
In the first phase the polynomial~$q$
is approximated by a trinomial with roots inside the interval, an application of using an adaptive approximation that resembles the function it 
approximates, while in the second phase, this trinomial is solved using adaptive approximation, where the approximation function
is constructed to ensure that the iterates are properly constrained. Here we concentrate on the second phase, as it provides a good
and uncomplicated example of adaptive approximation designed to satisfy specific requirements on the iterates.

The trinomial equation we need to solve is given by $f(x) \coloneqq ax^{n} - bx^{k} + c =0$, where $a,b,c >0$, $n \geq 3$, and $1 \leq k  \leq n-1$,
under the assumption that~$f$ has two positive roots~$\rone$ and~$\rtwo$ with $\rone < \rtwo$. 
The goal is to find an approximation to~$f$ that dominates it so that the approximation has roots in the interval~$[\rone,\rtwo]$. To facilitate
this, we use the transformation of variables~$z=x^{k}$, which transforms~$f$ into $F(z)=f(z^{1/k})=az^{n/k} - bz +c$. The function~$F$
is convex, so that using a linear approximation is not appropriate as it would be dominated by~$F$. On the other hand, $z^{-n/k}$ is also convex, 
so that it dominates its linear approximation. Consequently, the reciprocal of this linear approximation then approximates $z^{n/k}$ to first order, 
and it dominates~$z^{n/k}$, which is precisely what we need, as will soon become clear. This is the general idea that we now consider in more detail.

\noindent {\bf \underline{Method construction.}} 
The first order approximation of $z^{n/k}$ at $z=\bar{z}$ by~$R(z)=\alpha/(\beta - z)$ is obtained by setting $R(\bar{z})=\bar{z}^{n/k}$ and
$R^{\prime}(\bar{z})=(n/k)\bar{z}^{n/k-1}$. A straightforward calculation shows that $\alpha=(k/n)\bar{z}^{1+n/k}>0$ and $\beta=(1+k/n)\bar{z}>\bar{z}>0$.
Since $R(z)$ approximates $z^{n/k}$ to first order, $1/R(z)$, which is linear, approximates $z^{-n/k}$ to first order. 
Since $z^{-n/k}$ is convex, this means that
\bdis 
\dfrac{1}{R(z)} = \dfrac{\beta-z}{\alpha} \leq z^{-n/k} \Longrightarrow R(z) = \dfrac{\alpha}{\beta -z} \geq z^{n/k} \; . 
\edis
As a result, we have obtained, with $G(z)=aR(z)-bz+c$, that $G(z) \geq F(z)$.
If $F(\bar{z})=G(\bar{z}) < 0$, then the approximation necessarily has roots in the interval $[\rone^{k},\rtwo^{k}]$, 
since~$G(0)>0$ and $G(z) \rightarrow +\infty$ as $z \rightarrow \beta$. There are two such roots as they are the solution of 
a quadratic equation, and they become the next iterates. 
Figure~\ref{figure_trinom} shows the trinomial~$f(x)$ and its approximation $g(x)=G(x^{k})$, as well as $F(z)$ and its approximation $G(z)$.

\noindent {\bf \underline{Convergence properties.}}
A method based on the approximation~$G$, starting from a point with negative function value, iterates with the smaller or larger root of~$G$, 
to converge monotonically to $\rone^{k}$ or~$\rtwo^{k}$, respectively, from within the interval. It is a direct consequence of the domination of~$F$ by~$G$.
The order of convergence is quadratic (\cite{Melman_Pellet}).                              

A starting point~$z_{0}$ is most conveniently computed as the minimum argument of~$F$, obtained from 
\bdis
F^{\prime}(z)=\dfrac{n}{k} a z^{n/k-1} - b = 0 \Longrightarrow z_{0} = \lb \dfrac{kb}{na} \rb^{\frac{k}{n-k}} \hskip -0.25cm .
\edis
\vskip -0.3cm   
%
%
%
%
\begin{figure}[H]
\begin{center}
\raisebox{0ex}{\includegraphics[width=0.325\linewidth]{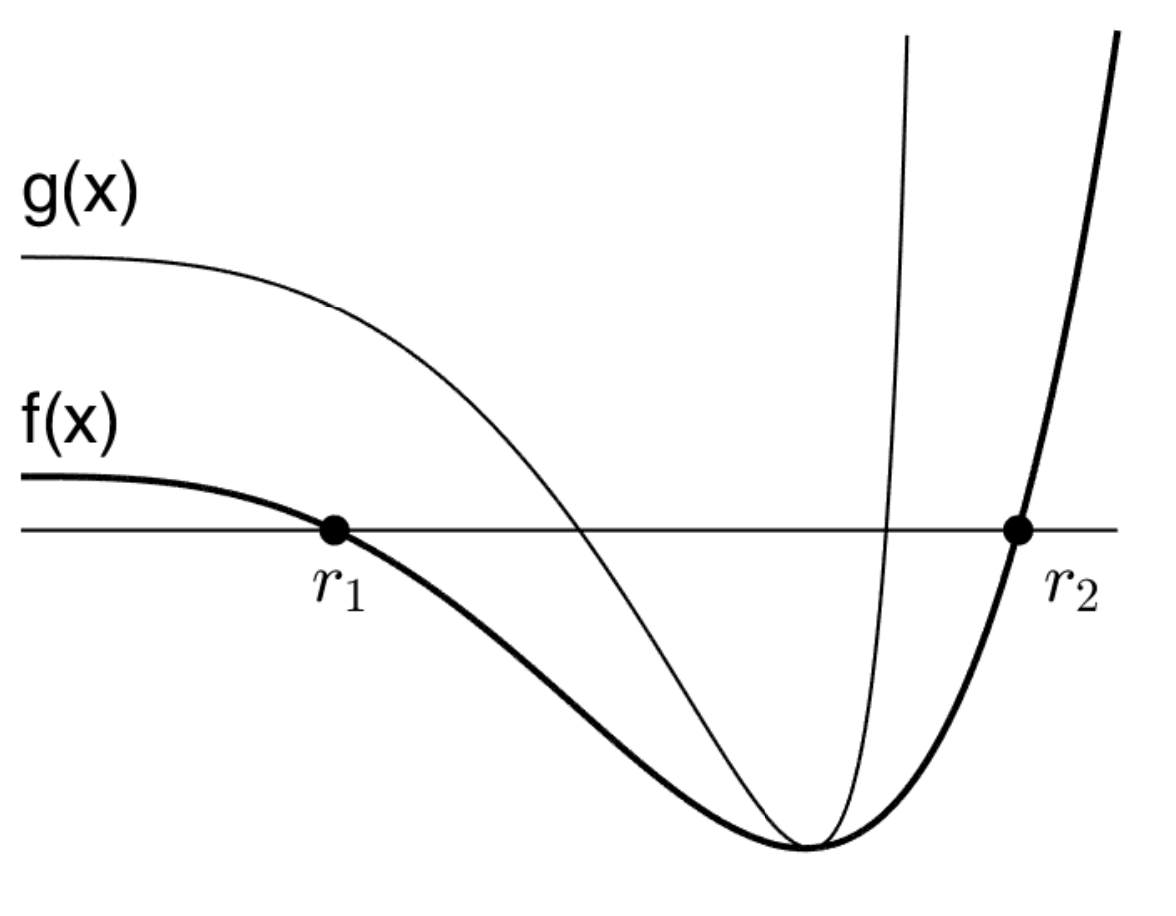}}
\hskip 2.0cm
\raisebox{1ex}{\includegraphics[width=0.325\linewidth]{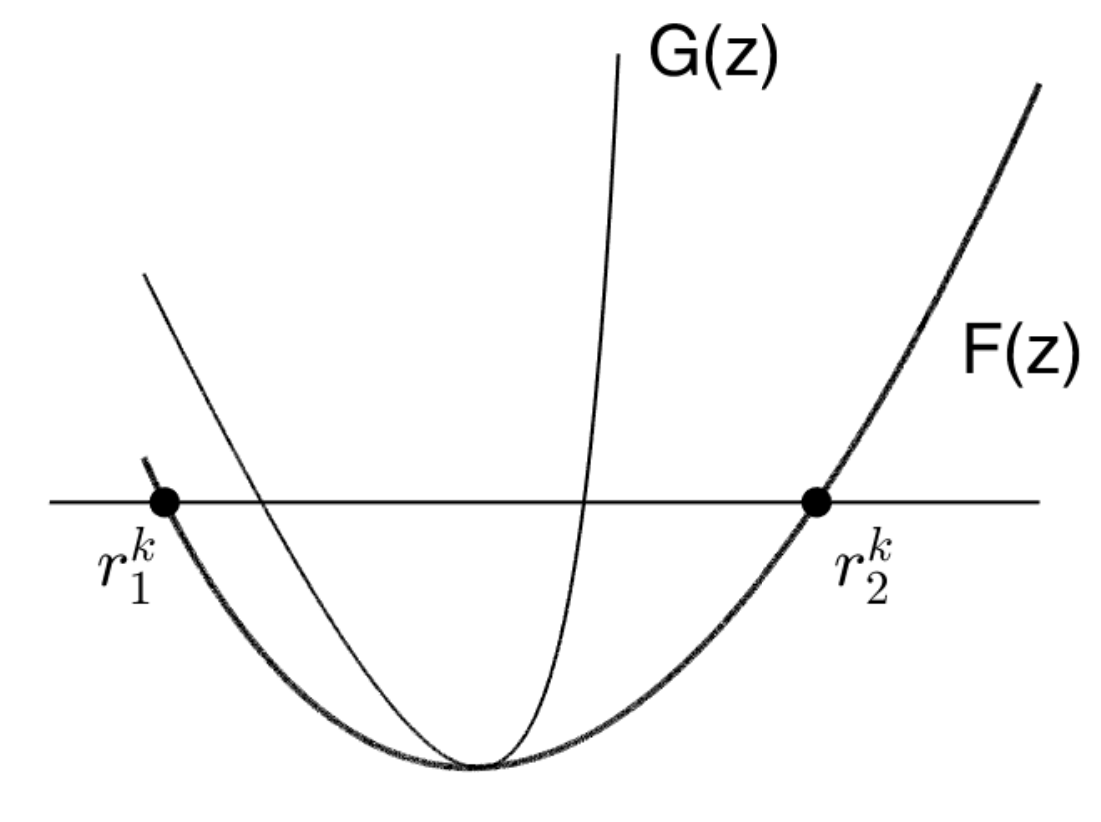}}
\caption{The functions~$f(x)$ and~$F(z)$ with their respective approximations~$g(x)$ and~$G(z)$.}
\label{figure_trinom}            
\end{center}
\end{figure}
\end{example}

\section{Conclusion}
Most numerical analysis textbooks treat the problem of solving a single nonlinear equation in a real variable by deriving and analyzing standard 
all-purpose methods that, unfortunately, do not work very well (or at all) for many problems. In contrast, we have demonstrated the importance 
of constructing taylor-made methods that take the given problem's characteristics into account by presenting four case studies
where standard methods are inefficient or cannot be used. In each case we constructed a specialized method using a few general tools that
can be helpful for other problem types as well.

\end{document}